\documentclass[12pt]{amsart}      
\usepackage{amsmath,amssymb,amsfonts} 
\usepackage{graphicx}                 

\theoremstyle{plain}
\newtheorem{theorem}{Theorem}[section]
\newtheorem{proposition}[theorem]{Proposition}
\newtheorem{corollary}[theorem]{Corollary}
\newtheorem{lemma}[theorem]{Lemma}
\newtheorem{conjecture}[theorem]{Conjecture}

\theoremstyle{definition}
\newtheorem{definition}[theorem]{Definition}
\newtheorem{problem}[theorem]{Problem}
\newtheorem{example}[theorem]{Example}
\newtheorem{A}[theorem]{Assumption}

\theoremstyle{remark}
\newtheorem{remark}[theorem]{Remark}

\def\b#1{\mathbf{#1}}
\def\f#1{\mathfrak{#1}}
\def\bb#1{\mathbb{#1}}
\def\c#1{\mathcal{#1}}
\def\getsfrom{\, \text{\rotatebox[origin=c]{180}{$\mapsto$} } }

\DeclareMathOperator{\Hilb}{Hilb}

\DeclareMathOperator{\Nash}{Nash}
\DeclareMathOperator{\fNash}{fNash}
\DeclareMathOperator{\Spec}{Spec}
\DeclareMathOperator{\ord}{ord}
\DeclareMathOperator{\initial}{in}
\DeclareMathOperator{\Hom}{Hom}

\date{}
\title{Flag higher Nash blowups}
\author{Takehiko Yasuda}
\thanks{Financially supported by the Japan Society for the Promotion of Science.}
\address{Research Institute for Mathematical Sciences, Kyoto University, Kyoto, 606-8502, Japan}
\email{takehiko@kurims.kyoto-u.ac.jp}
\begin{document}

\maketitle

\begin{abstract}
In his previous paper \cite{highernash}, the author has defined
a higher version of the Nash blowup and considered it a possible
candidate for the one-step resolution.
In this paper, we will introduce another higher version of the Nash blowup
and prove that it is compatible with products and smooth morphisms.
We will also prove that the product of curves can be desingularized
via both versions. 
\end{abstract}

\section{Introduction}

Let $X$ be a variety over an algebraically closed field $k$
of characteristic 0, and  $X_{sm}$ its smooth locus.
For $x \in X$ and a non-negative integer $n$, put $x^{(n)}:= \Spec \c O_{X} /\f m_{x}^{n+1}$.
The author \cite{highernash} has defined the {\em $n$-th Nash blowup} of $X$, 
denoted $\Nash _{n} (X)$,
 to be the closure of the set 
$ \{ x^{(n)} | x \in X_{sm} \} $ in the Hilbert scheme $\Hilb (X)$.
In this paper, we call it the {\em simple $n$-th Nash blowup}, distinguished
from what we will introduce below.
Unfortunately, as we will see (Example \ref{expl-not-isom}),
 it is not generally compatible with either products or smooth morphisms.
This would be explained as follows: Let $Y$ be another
variety,  and  let $U \subset X$ and $V \subset Y$
be clusters (that is, 0-dimensional closed subschemes)
with $U \in \Nash _{n} (X)$ and $V \in \Nash _{n}(Y)$.
Then $U \times V$ seems the only natural cluster in $X \times Y$ constructed
from $U$ and $V$.
But it is not of the expected length $\binom{\dim X + \dim Y +n}{n}$, 
so $U \times V \notin \Nash _{n}(X \times Y)$.
 In particular,
for $x \in X_{sm}$ and $y \in Y_{sm}$, 
$x ^{(n)} \times y ^{(n)} \ne (x,y)^{(n)}$.

For instance, the ideal $(u,v)^{n+1} \subset k[u,v]$
 is not identical to either $(u^{n+1},v^{n+1})$ or $(u^{n+1})(v^{n+1})$,
but to $\bigcap _{i=1}^{n+1} (u^{i})+(v^{n+2-i})$,
$\sum _{i=0}^{n+1} (u^{i}) \cap (v^{n+1-i})$
and $\sum _{i=0}^{n+1} (u^{i}) (v^{n+1-i})$.
This observation suggests considering a collection
 $(x ^{(0)},\dots, x^{(n)})$
rather than a single $x^{(n)}$.
We now define the {\em flag $n$-th Nash blowup}, denoted $\fNash _{n} (X)$,
to be the closure of the set 
$ \{ (x^{(0)},\dots,x^{(n)}) | x \in X_{sm}  \} $
in $ ( \Hilb (X) )^{n+1}$. 
This is also a higher version of the classical Nash blowup. Namely 
$\fNash _{1} (X)$ is canonically isomorphic to the classical Nash blowup.
By definition, every point of $\fNash _{n} (X)$
is the collection  of $n+1$ clusters in $X$.

For $U _{*} =(U_{0},\dots,U_{n})\in \fNash _{n} (X)$ and $V _{*} =(V_{0},\dots,V_{n})
\in \fNash _{n} (Y)$
and for $0 \le i \le n$, put
\[
 \{U_{*},V_{*}\}_{i}:= \bigcup_{j= 0}^{i}  U_{j}  \times V_{i -j} \subset X \times Y ,
\]
and $\{ U_{*},V_{*} \}_{*} := ( \{ U_{*},V_{*} \}_{0}, \dots, \{ U_{*},V_{*} \}_{n} ) $.
We will see $\{ U_{*} ,V_{*}\}_{*} \in \fNash _{n} (X \times Y) $.
The main theorem of this paper is the following.

\begin{theorem}\label{main-thm}
Let $X$ and $Y$ be varieties.
Then there is a canonical isomorphism, 
\begin{align*}
 \fNash _{n} (X) \times \fNash _{n} (Y) & \cong \fNash _{n} (X \times Y) \\
 (U _{*} ,V _{*}) & \mapsto \{ U _{*} , V _{*} \}_{*}.
\end{align*}
\end{theorem}

Using this, we will prove also that 
the flag higher Nash blowup is compatible with smooth morphisms.
In general, it is expected that a good resolution is compatible
with smooth morphisms. Resolution with this property  was first constructed
by Villamayor \cite{VillamayorConstructiveness, VillamayorPatching}
(see also \cite{BierstoneMilmanCanonicalResolution, EncinasVillamayorActa}).
This property enables us to even construct resolution of Artin stacks.

The second half of the paper will be devoted to the study of 
simple and flag higher Nash blowups of 
curves and products of curves. 
Consider a formal irreducible curve $\b X = \Spec A$ (that is, $A$
is the complete Noetherian local domain of dimension 1 with coefficient
field $k$).  We can also define  simple and flag higher Nash blowups of 
such ``formal varieties''.
Fix an embedding $A \hookrightarrow k[[x]]$ so that 
$k[[x]]$ is the integral closure of $A$.
Then
the {\em associated numerical monoid} of $\b X$ (of $A$) is
\[
S := \{ s \in \bb Z_{\ge 0} | \exists f \in A, \ \ord f = s  \} = \{0= s_{0} < s_{1}< \cdots  \} \subset \bb Z_{\ge 0}.
\]
Then a theorem
in \cite{highernash} says that 
$\Nash _{n} (\b X)$ is  normal if and only if $s_{n+1}-1 \in S$.
(Notice that the indices of $s_{i}$ are shifted by one from those in \cite{highernash}.)
For the flag higher Nash blowup, we will obtain a similar result;
 $\fNash _{n} (\b X)$ is normal if and only if 
$s_{m} -1 \in S$ for some $m \le n+1$.

Next consider formal curves $\b X_{i}$, $1 \le i\le l$.
Let 
\[
S _{i}= \{0= s_{i,0} < s_{i,1}< \cdots  \} \subset \bb Z_{\ge 0}
\]
be their associated numerical monoids respectively. 
It follows from Theorem \ref{main-thm} that
$\fNash _{n} (\hat \prod _{i} \b X_{i})$
is regular if and only if $\forall i$, $\exists m_{i} \le n+1$,
$s_{i,m_{i}}-1 \in S_{i}$.
Quite strangely, when $l \ge 2$, 
$\Nash _{n} (\hat \prod _{i} \b X_{i})$
is regular exactly when so is $\fNash _{n} (\hat \prod _{i} \b X_{i})$,
though its proof is more involved.

Let $X$ be a variety whose analytic branches all have
the singularity type of the product of curves,
$C \subset X$ its conductor subscheme,
$Z \in \Nash _{n}(X)$ and $Z_{*} \in \fNash _{n}(X)$. 
Then it is a corollary of the preceding results
that if $Z \not \subset C$ (resp.\ $Z_{n}\not \subset C$),
then $\Nash _{n} (X)$ (resp.\ $\fNash _{n} (X)$) is smooth
at $Z$ (resp.\ $Z_{*}$).  

The simple $(n+1)$-th Nash blowup is not generally isomorphic 
to the $n$-th even if the $n$-th is smooth. So the simple $n$-th Nash blowup is not generally
isomorphic to the $n$-times iteration of classical Nash blowups.
On the other hand, at least for curves and products of curves,
which are the only cases computed up to now, 
the flag $(n+1)$-th Nash blowup is  isomorphic 
to the $n$-th if the $n$-th is smooth. Moreover the flag higher Nash blowup 
and the iteration of classical Nash blowups are both 
compatible with products and smooth morphisms.
However, we compute various blowups
of $\Spec k[[x^{5},x^{7}]]$ in \S \ref{5-7-example}, which shows that
the flag $n$-th Nash blowups and the $n$-times iteration
of the classical blowups are not generally isomorphic to each other.

\subsection{Convention}

Throughout the paper, we work over
an algebraically closed base field $k$ of characteristic zero. 
The product of $k$-schemes, denoted $\times$, means 
the fiber product over $k$, and the tensor product of $k$-algebras, 
denoted $\otimes $, means the tensor product over $k$.
A {\em point} of a scheme means a $k$-point.
A {\em variety} means a separated integral scheme of finite type over $k$. 
A {\em cluster} means a zero-dimensional closed subscheme.
For a variety $X $, we denote by $\Hilb (X)$ the Hilbert scheme of $X$ of clusters.
For a cluster $Z \subset X$, we write $Z \in \Hilb (X)$ for the corresponding
point. 
Similarly, if $Z_{*} := (Z_{0},\dots,Z_{n})$ is a collection of clusters,
 then we write
 $Z_* \in \Hilb (X) ^{n+1}$.

\section{Simple and flag higher Nash blowups}

In this section, we define the simple and flag higher Nash blowups
and show their basic properties. For more details of simple
higher Nash blowup, we refer the reader to \cite{highernash}.

Let $X$ be a variety of dimension $d$. For a point $x \in X$ with 
defining ideal $\f m \subset \c O_{X}$,
its $n$-th infinitesimal neighborhood, denoted $x^{(n)}$, is
the cluster defined by $\f m ^ {n+1}$.
We denote by $X_{sm}$ the smooth locus of $X$.

\begin{definition}\label{def-simple}
The {\em simple $n$-th Nash blowup}, denoted $\Nash _n (X)$, is the closure
of the set $ \{   x ^{(n)} | x \in X_{sm}  \} $ in $\Hilb (X)$ (endowed with 
the reduced scheme structure).
\end{definition}

\begin{remark}
In \cite{highernash}, $\Nash _{n} (X)$ was defined to be the closure
of $\{ (x,x^{(n)}) | x \in X_{sm} \} \subset X \times \Hilb (X)$,
which is however canonically isomorphic to $\Nash _{n}(X)$ in Definition 
\ref{def-simple} in characteristic zero.
\end{remark}

Every point $Z \in \Nash _n (X)$ is a cluster whose coordinate ring
is a $k$-vector space of dimension $\binom{n+ d}{ n }$.
Moreover as a set, $Z$ is a point, say $x \in X$. We have $Z \subset x^{(n)}$.
But the equality does not generally holds unless $x \in X_{sm}$. 
There is a natural morphism
\[
  \Nash _n (X) \to X , \, Z \mapsto x ,
\]
which is the morphism forgetting the scheme structure of $Z$.
The morphism is projective and an isomorphism over $X_{sm}$.

\begin{definition}
 The {\em flag $n$-th Nash blowup}, denoted $\fNash_n (X)$, 
is  the closure of the set $
 \{  ( x ^{(0)} , x ^{(1)} , \dots , x^{(n)} )|  x \in X_{sm}\} $
in $\Hilb (X) ^{n+1}$.
\end{definition}

Every point of  $\fNash _n (X)$ is a collection $ Z_* =(Z_{0},Z_{1},\dots,Z_{n})$ such that 
 $Z_i$ is a cluster of length $ \binom{  i + d }{i } $ and
$Z_ i \subset Z_{i+1}$, $ 0 \le i \le n -1$.
In particular, $Z_0$ is a reduced point and we have a map
\[
  \fNash _n (X) \to X, \, Z_* \mapsto Z_0.
\]
Again this is  projective and an isomorphism over $X_{sm}$.
 Moreover for $ m \le n $, there are  natural projections
\begin{gather*}
 \fNash _n (X) \to \fNash_m (X), ( Z_0 , \dots, Z_n  ) \mapsto
(Z_0 , \dots, Z_m ) , \text{ and}\\
 \fNash _n (X) \to \Nash_m (X), ( Z_0 , \dots, Z_n  ) \mapsto
 Z_m ,
\end{gather*}
which are projective and birational.
It is easy to see that $\fNash _n (X)$ is isomorphic to the irreducible component of 
\[
 \Nash_n (X) \times_X \Nash_{n-1} (X) \times_X \cdots \times_X 
\Nash_1 (X) \times_X \Nash_0 (X)
\]
which dominates $X$. 
Indeed, they are the same subscheme of $ ( \Hilb (X) )^{n+1} $.

There is a slightly different construction
of higher Nash blowups via relative Hilbert scheme.
Let $\Delta \subset X \times X $ be the diagonal
and $\Delta^{(n)}$ its $n$-th infinitesimal neighborhood.
Namely if $I_{\Delta}$ is the defining ideal sheaf of 
$\Delta$, then $\Delta^{(n)}$ is the closed subscheme of $X \times X$ defined by
$I ^{n+1}_{\Delta}$.
We think of $\Delta^{(n)}$ as an $X$-scheme via the first projection.
If $ \Hilb _{ \binom{n+d}{n}} (\Delta ^{(n)} / X)$ is the relative Hilbert scheme
of clusters of length $\binom{n+d}{n}$,  then 
$\Nash_{n} (X)$ is canonically isomorphic to the irreducible
component of  $ \Hilb _{ \binom{n+d}{n}} (\Delta ^{(n)} / X)$
dominating $X$.

We adopt this construction as the definition for formal varieties.
Let $A$ be a complete local Noetherian reduced ring with coefficient field $k$
and put $\b X:= \Spec A$.
Let 
$\b  {\Delta} ^{(n)} \subset \b X \hat \times \b X := \Spec A \hat \otimes A$
be 
the $n$-th infinitesimal neighborhood of the diagonal.
We regard this as an $\b X$-scheme via the first projection.
We make the following assumption:
\begin{A}\label{Assumption}
$\b X$ has pure dimension $d$ and there exists an open dense subscheme of $ X$
over which $\b  \Delta^{(n)}$ is flat and finite of degree
$\binom{n+d}{d}$.
\end{A}

The author does not know whether the second condition in Assumption \ref{Assumption} always holds under the first.
But it holds at least in the following case;
$\b  X'$ is the completion of a variety $X$ at a point,
$\{ \b  X_{i} \}_{i \in I}$ a subcollection of irreducible components of $\b X$
and $\b X = \cup_{i \in I} \b X_{i}$. 
In particular, in the case where $\dim \b X =1$.

\begin{definition}
We define the {\em simple (resp.\ flag) $n$-th Nash blowup} of $\b X$,
denoted $\Nash _{n} (\b X) $ (resp.\ $\fNash _{n} (\b  X)$), to be the union of those irreducible components of 
 $\Hilb _{\binom{d+n}{n}}(\b {\Delta} ^{(n)} /  \b X)$
\[
  (\text{resp.\ } \Hilb _{\binom{d+n}{n}}(\b  {\Delta} ^{(n)} / \b  X) \times _{ \b X} \Hilb _{\binom{d+n-1}{d}}(\b  {\Delta} ^{(n-1)} /\b   X) \times _{\b  X} \cdots \times _{\b  X} \Hilb_{1} (\b \Delta / \b X) )
\]
which dominate irreducible components of $\b  X$. 
\end{definition}

If $X$ is a variety and $\b X$ is its completion at a point,
and if $\Delta \subset X \times X$ and $\b \Delta \subset \b X \hat \times \b X$
are the diagonals respectively, 
then 
$\Delta ^{(n)} \times _{X} \b X \cong \b \Delta ^{(n)}$,
so $\Hilb (\Delta /X) \times_{X} \b X \cong \Hilb (\b \Delta/ \b X)$.
Hence if $\b X_{i}$, $i \in \Lambda$, are the irreducible components of $\b X$, 
then 
\begin{gather*}
\Nash _{n} (X) \times _{X} \b X \cong \Nash _{n} (\b X) = \bigcup_{i \in \Lambda} \Nash _{n} (\b  X_{i}) \\
\fNash _{n} (X) \times _{X} \b X \cong \fNash _{n} (\b  X) = \bigcup_{i \in \Lambda} \fNash _{n} (\b  X_{i}) 
\end{gather*}

\section{Compatibility with products}

Let $X $ and $Y$ be varieties of dimension $d$ and $e$ respectively, 
and $U_{*} \in \fNash_{n}(X), V_{*}\in \fNash_{n}(Y)$. 
We follow  the convention that $U_{-1} := \emptyset$ and $V_{-1} := \emptyset$.
For closed subschemes $A,B \subset C$ defined by the ideal sheaves
$I,J \subset \c O_{C}$ respectively, define $A \vee B \subset C$
to be the closed subscheme defined by $IJ$.

\begin{proposition}\label{prop-length-product}
For each $ 0 \le i \le n $,
we have the following identification of clusters of $X \times Y$,
\[
 \bigcup_{j= 0}^{i}  U_{j}  \times V_{i -j} = 
 \bigcap _{j=-1}^{i}  ( U_{j}  \times Y) \cup  (X \times V_{i-1 -j} )
 =\bigcap _{j=-1}^{i}  ( U_{j}  \times Y) \vee  (X \times V_{i-1 -j} ).
 \]
Moreover this cluster is of length $\binom{d+e+i}{i}$.
When $U_{*} =(x^{(0)},\dots,x^{(n)})$ and 
$V=(y^{(0)},\dots,y^{(n)})$ with $x \in X_{sm}$
and $y \in Y_{sm}$, this cluster is identical to
$(x,y)^{(i)}$.
\end{proposition}

\begin{proof}
We may suppose that 
$X$ and $Y$ are affine, say
 $ X= \Spec R$, $Y= \Spec S$.
Let $\f m _{i+1} \subset R$ and $\f n_{i+1} \subset S$
be the defining ideals of $U_{i}$ and $V_{i}$ ($i=0,1,\dots, n$) respectively.
Put $\tilde{\f m}_{i}:= \f m_{i} (R \otimes S)$ and $\tilde{\f n}_{i}:= \f n_{i}( R \otimes S)$.
Then the left side of the equation in the proposition
is defined by the ideal
 $\bigcap _{ j =1} ^{i+1} \tilde{\f m}_{j} + \tilde {\f n}_{i+2-j} $,
 the middle by  $\sum _{j=0}^{i+1}   \tilde{\f m}_{j} \cap \tilde {\f n}_{i+1-j} $,
 and the right by $\sum _{j=0}^{i+1}   \tilde{\f m}_{j}  \tilde {\f n}_{i+1-j} $.

Take bases $ A \subset R$
(resp.\ $B \subset S$) of $R$
(resp.\ $S$) as $k$-vector spaces such that
for each $i$, $\{ a \in A| a \in \f m _{i}  \}$
(resp.\ $\{ b \in B| b \in \f n _{i}  \}$)
is a basis of $\f m _{i}$ (resp.\ $\f n_{i}$).
Then $\{ a \otimes b | a \in A , \ b \in B\}$
is a basis of $R \otimes S$.
For $f \in S$ and $g \in R$, define 
\begin{gather*}
\ord  f := \max \{ i| f \in \f m_{i}\} \in \{0,1,\dots, n+1 \} , \\
 \ord  g := \max \{ i| g \in \f n_{i}\} \in \{0,1,\dots, n+1 \} , \text{ and}\\ 
  \ord f\otimes g := \ord f + \ord g \in \{0,1,\dots,2 (n+1) \}  .
\end{gather*}

For $h = \sum_{i=1}^{l}  c _{i}a _{i} \otimes b _{i}$
($c_{i} \in k \setminus \{0\}$, $a_{i} \in A$, $ b_{i} \in B$),
$\ord h := \min_{i} (\ord a_{i} \otimes b_{i}  )$.
We claim that
\[
\{ h | \ord h \ge i+1 \}=
\sum _{j=0}^{i+1}   \tilde{\f m}_{j}  \tilde {\f n}_{i+1-j} 
=\sum _{j=0}^{i+1}   \tilde{\f m}_{j} \cap \tilde {\f n}_{i+1-j} 
=\bigcap _{ j =1} ^{i+1} \tilde{\f m}_{j} + \tilde {\f n}_{i+2-j} .
\]
It is easy to see that 
\[
\{ h | \ord h \ge i+1 \}
\subset
\sum _{j=0}^{i+1}   \tilde{\f m}_{j}  \tilde {\f n}_{i+1-j} 
\subset \sum _{j=0}^{i+1}   \tilde{\f m}_{j} \cap \tilde {\f n}_{i+1-j} 
\subset \bigcap _{ j =1} ^{i+1} \tilde{\f m}_{j} + \tilde {\f n}_{i+2-j} .
\]
Since $ \tilde {\f m} _{j} + \tilde{\f n}_{i +2 -j} = 
 \langle a\otimes b | a \in A \cap \f m _{j} \text{ or } b \in B \cap \f n_{i+2-j} \rangle$,
 we have
\begin{align*}
&\bigcap _{ j =1} ^{i+1} \tilde{\f m}_{j} + \tilde {\f n}_{i+2-j}  \\
&= \langle a\otimes b | \forall j,\ a \in A \cap \f m _{j} \text{ or } b \in B \cap \f n_{i+2-j} \rangle \\
& = \langle a\otimes b | \forall j,\ \ord a \ge j \text{ or } \ord b \ge i+2-j \rangle \\
&= \langle a\otimes b |\ord b \ge i+2-(\ord a +1) \rangle \\
& = \langle a\otimes b |  \ord a \otimes b \ge i+ 1 \rangle \\
&= \{ h | \ord h \ge i+1 \} .
\end{align*}
This proves our claim and the first assertion of the proposition.

The length of our cluster is equal to 
$\sharp \{a \otimes b |a \in A, \ b \in B, \ \ord a \otimes b \le i \}$,
which depends only on $d$, $e$ and $i$.
So the second assertion follows from the last.

To show the last assertion,
we may assume that $R=k[x_{1},\dots,x_{d}]$
and $S=k[y_{1},\dots,y_{e}]$,
and that $\f m _{i} =(x_{1},\dots,x_{d})^{i}$
and $\f n_{i} = (y_{1},\dots,y_{e})^{i}$.
Then for $h \in R \otimes S =k[x_{1},\dots,x_{d},y_{1},\dots,y_{e}]$,
the order defined above, $\ord h$, is equal to the usual order of
polynomial. 
So $\{ h |\ord h \ge i+1 \} =(x_{1},\dots,x_{d},y_{1},\dots,y_{e})^{i+1}$,
which completes the proof. 
\end{proof}

\begin{definition}
We denote the cluster in Proposition \ref{prop-length-product}
by $\{ U_{*},V_{*} \}_{i}$
and put $ \{ U_{*},V_{*} \}_{*}:=(\{ U_{*},V_{*} \}_{1},\dots,\{ U_{*},V_{*} \}_{n}) $.
\end{definition}

\begin{lemma}\label{lem-product-well-def}
There exists a morphism 
\[
\mu : \fNash _n (X) \times \fNash _n (Y) \to \fNash _n ( X \times Y)
\]
that takes $ (U _* , V _*) $ to $ \{ U_{*}, V _{*} \} _* $.
Moreover it is birational and surjective.
\end{lemma}

\begin{proof}
Let $ \c U _*  := ( \c U_1 , \dots, \c U _{n+1} )$, 
$\c U_{i} \subset \fNash _{n}(X) \times X$
 (resp.\ $ \c V _*  := ( \c V_1 , \dots, \c V _{n+1} )$, $\c V_{i} \subset \fNash _{n}(Y)\times Y$) be the universal collection of clusters over $\fNash _n (X)$
(resp.\ $\fNash _n (Y)$).
Set 
\[
 \{\c U _{*},\c V_{*}\} _i := \bigcap _{j=-1}^{i}  (\c U_{j} \times \fNash _{n}(Y)\times Y) \vee (\fNash _{n}(X) \times X \times \c V_{i-1-j}).
\]
Then the fiber of the projection
$ \{\c U _{*},\c V_{*}\} _i  \to \fNash _{n}(X) \times \fNash _{n} (Y)   $ 
over $(U_{*},V_{*})$ is  
$ \{U_{*},V_{*}  \}_{i} $. 
So from Proposition \ref{prop-length-product},
 $\{\c U _{*},\c V_{*}\} _i$ is flat over $\fNash _{n} (X) \times \fNash _{n} (Y)$ and 
 generically the family of $z^{(i)}$, $z \in X_{sm}\times Y_{sm}$. 
 From the universality, there exists a morphism 
$ \fNash _n (X) \times \fNash _n (Y) \to \fNash _n ( X \times Y) $
corresponding to $\{\c U _{*},\c V_{*}\} _*$, which takes 
 $ (U _* , V _*) $ to $ \{ U _{*}, V_{*}\} _* $.
It is an isomorphism over $X_{sm} \times Y_{sm}$.
 Since $  \fNash _n (X) \times \fNash _n (Y) $
is proper over $X \times Y$, the morphism is surjective.
\end{proof}

Let $ p :X \times Y \to X $ and $q : X \times Y \to Y$ be the projections.
For $Z_* = ( Z_0, \dots, Z_n ) \in \fNash _n (X \times Y)$, we denote by 
$p (Z_i )$ (resp.\  $q (Z_i)$) the scheme-theoretic image of $Z_i$
by $p$ (resp.\ $q$), and set $p (Z_*) := ( p (Z_0) , \dots, p(Z_n)  )$
and $q (Z_*) := ( q (Z_0) , \dots, q(Z_n)  )$.

\begin{theorem}\label{thm-product}
We have a canonical isomorphism
\begin{align*}
\fNash _n (X) \times \fNash _n (Y) &\cong \fNash _n ( X \times  Y) \\
       ( U _ * , V_ *  ) & \mapsto \{ U_{*} ,V _{*}\} _* \\
          ( p ( Z_* ), q (Z_ *)  ) &  \getsfrom  Z_* .
\end{align*}
\end{theorem}

\begin{proof}
Since  $\mu $ in Lemma \ref{lem-product-well-def} is surjective, any $Z_*  \in \fNash _n (X \times Y)$ is
 $\{ U_{*} ,V _{*}\} _*  $ for some $(U _* , V _* ) \in \fNash _n (X) \times \fNash _n (Y)$.
Then it is clear that 
$p(Z_{*})=U_{*}$ and $q(Z_{*})=V_{*}$.
Thus the map $\fNash _{n } (X \times Y) \to \fNash _{n } (X) \times \fNash _{n } (Y)$,
$Z_{*} \mapsto (p(Z_{*}),q(Z_{*}) )$ is the inverse of $ \mu$.
We will show this map is a morphism of schemes.

We have a natural morphism
$ \fNash _{n} (X \times Y) \to X \times Y \to Y $,
which induces a closed embedding 
\[
\fNash _{n} (X \times Y) \times X \hookrightarrow \fNash _{n} (X \times Y) \times X \times Y .
\]
We denote the image of this embedding by $W$.
Let $\c Z_{*} \subset \fNash _{n} (X \times Y) \times X \times Y$ be 
the universal collection of clusters.
Then the scheme-theoretic intersection $\c Z_{i} \cap W$
is a family of clusters in $X$ over $\fNash _{n} (X \times Y)$.
If $F$ is the fiber of the projection 
$\fNash _{n} (X \times Y) \times X \times Y \to \fNash _{n} (X \times Y)$
over $Z_{*}$,
then the fiber of $\c Z _{i}\cap W \to \fNash_{n} (X \times Y)$
over $Z _{*}$ is
\[
 \c Z_{i} \cap W   \cap F = (\c Z_{i} \cap F  ) \cap (W \cap F) = Z_{i} \cap  (X \times q(Z_{0})) = p (Z_{i}) .
\]
So $\c Z_{i} \cap W$ is a flat family, generically of 
$x^{(i)}$, $x \in X_{sm}$.
From the universality, there exists the corresponding morphism
$\fNash _{n} (X \times Y) \to \fNash _{n} (X)$.
Similarly we obtain $ \fNash _{n} (X \times Y) \to \fNash _{n} (Y)$,
and $ \fNash _{n} (X \times Y) \to \fNash _{n} (X) \times \fNash _{n} (Y)$.
The last morphism is clearly the inverse of $\mu$.
\end{proof}

It is straightforward to generalize the  theorem above
to the product of an arbitrary number of varieties.
Let $X_{i}$, $i=1,\dots, m$, be varieties.
For $U _{i,*} \in \fNash _{n} (X_{i}) $, $i=1,\dots,m$,
we set 
\begin{align*}
 \{ U _{1, *} ,\dots,U _{m , *}\} _{j}
 & := 
 \{ \{ \cdots \{\{ U_{1,*},U_{2,*} \}_{*},U_{3,*}\}_{*}, \cdots \}_{*},U_{m,*}\} _{j}\\
& = \bigcup _{\sum _{ i} j _{i} =j} \prod _{i} U_{i, j_{i}} .
\end{align*}
Let $p_{l} : \prod_{i} X_{i} \to  X_{l}$
be the $l$-th projection.

\begin{theorem}
We have a canonical isomorphism
\begin{align*}
\prod _{i}\fNash _n (X_{i})  &\cong \fNash _n ( \prod _{i}X _{i}) \\
       ( U _ {1,*} ,\dots, U _{m,*}  ) & \mapsto \{ U_{1,*} , \dots ,U _{m,*}\} _* \\
          ( p_{1} ( Z_* ),\dots, p_{m} (Z_ *)  ) &  \getsfrom  Z_* .
\end{align*}
\end{theorem}

\section{Compatibility with smooth morphisms}

If $f:Y \to X$ is an etale morphism of varieties, then there is a
natural isomorphism $\Delta^{(n)}_{Y} \cong \Delta ^{(n)}_{X} \times _{X} Y$,
where $\Delta_{Y} \subset Y \times Y$ and $\Delta_{X} \subset X \times X$
are the diagonals. So both simple and flag higher Nash blowups
are compatible with etale morphism, that is, there are canonical isomorphisms,
\begin{gather*}
 \Nash _{n} (Y) \cong \Nash _{n} (X) \times _{X} Y \\
  \fNash _{n} (Y) \cong \fNash _{n} (X) \times _{X} Y .
\end{gather*}
The composite $\fNash _{n} (Y) \to \fNash _{n} (X) \times _{X } Y \to
\fNash _{n} (X)$ takes $Z_{*}$ to $f(Z_{*})$.

Next if $f :Y \to X$ is a smooth morphism,
then there is an open covering $Y = \bigcup Y_{e}$ such that
for each $e$, $f |_{Y_{e}}$ factors as 
\[
 Y_{e} \xrightarrow{h} X \times \bb A^{c} \xrightarrow{p} X ,
\]
where $h$ is etale and $p$ is the projection. 
The morphisms 
\begin{gather*}
\fNash _{n}(Y_{e}) \to \fNash_{n} (X \times \bb A^{c}), \, 
Z _{*} \mapsto h (Z_{*}), \text{ and} \\
\fNash_{n} (X \times \bb A^{c}) \to \fNash _{n}(X), \, 
W _{*} \mapsto p (W_{*})
\end{gather*} 
are well-defined, and so is
\[
 \fNash _{n } (Y_{e}) \to \fNash _{n} (X) , \, Z_{*} \mapsto f (Z_{*}). 
\]
Gluing them yeild the morphism 
\[
 \fNash _{n } (Y) \to \fNash _{n} (X) , \, Z_{*} \mapsto f (Z_{*}).
\]

\begin{corollary}
Let $f : Y \to X$ be a smooth morphism of varieties. 
Then there is a canonical isomorphism 
\begin{align*}
 \fNash _n (Y) &\cong \fNash _n (X) \times _X Y  \\
 Z_{* } & \mapsto f(Z_{*}).
\end{align*}
\end{corollary}

\begin{proof}
Let the $Y_{e} $ be as above. Then
\begin{align*}
 \fNash _{n} (Y _{e}) & \cong \fNash _{n} (X \times \bb A^{c}) 
 \times _{X \times \bb A^{c} } Y_{e} \\
& \cong  (\fNash _{n} (X) \times \bb A^{c} ) \times _{X \times \bb A^{c} } Y _{e} \\
& \cong \fNash _{n} (X) \times _{X} Y_{e}.
\end{align*}
Since the isomorphisms $\fNash _{n} (Y_{e}) \cong \fNash _{n} (X) \times _{X} Y_{e}$
are canonical, we can glue them and obtain $ \fNash _n (Y) \cong \fNash _n (X) \times _X Y $.
\end{proof}

\section{Curves and products of curves}

\subsection{Curves}

Let $\b X = \Spec A$ be a formal irreducible curve. 
Namely $A$ is a complete Noetherian local domain with coeffficient field $k$.
Fix an embedding $A \hookrightarrow k[[x]]$ so that
 $k[[x]]$ is the integral closure of $A$. 
We define the {\em associated numerical monoid}
of $\b X$ (and $A$)
to be $ S := \{ s \in \bb Z_{\ge 0} | \exists f \in A, \, \ord f = s \} $ and  write
\[
 S = \{ 0= s_{0} < s_1 < s_2 < \cdots  \} .
\]
(Caution: The indices of $s_{i}$ differ by 1 from those in \cite{highernash}.)

\begin{theorem}[\cite{highernash}]\label{thm-formal-curve}
$\Nash _n (\b X) $ is normal if and only if $s_{n+1} -1 \in S $. 
\end{theorem}

\begin{proof}[Sketch of the proof]
Let $ \iota :  k[[x ]] \to k[[y]]$ be an isomorphism defined by $x \mapsto -y$.
The composite $A \hookrightarrow k[[x]] \xrightarrow{\iota} k[[y]]$ corresponds to a morphism 
$ \nu : \b Y := \Spec k[[y]] \to \b X$, which is the normalization of $\b X$. 
For each $n$, there exists a natural factorization of $\nu$ as follows; 
$\b  Y \xrightarrow{\phi_n} \Nash _n (\b X) \to \b X$. 
Then
$\phi_n$ corresponds to a family of clusters over $\b Y$, say $\c Z_{n} \subset \b X \hat \times \b Y := \Spec A[[y]]$.

Let $I \subset A[[y]]$ be the prime ideal defining the graph $\Gamma \subset\b  X \hat \times \b Y $ of $\nu$.
Then the defininig ideal of $\c Z_{n}$ is the $(n+1)$-th symbolic power $I^{(n+1)}$ of $I$. 
Let $\epsilon : \Spec k[y]/(y^2) \to \Spec k[[y]]$ be the natural morphism and set
$\c Z_{n, \epsilon } := \c Z_{n} \times _{\b Y} \Spec k[y]/(y^2) \subset \Spec A[y]/(y^2)$ and 
$Z_{n} :=  \c Z_{n} \times _{\b Y} \Spec k $. 
Put  
$\f a _{n+1} := I^{(n+1)} (A[y]/(y)) \subset A$, which is the defining ideal
of $Z_{n}$.
Then $\Nash _n (\b X)$ is normal if and only if
$\phi_{n}\circ \epsilon$ is a nonzero tangent vector 
if and only if
 $\c Z_{n,\epsilon}$ is a non-trivial family. 
We can construct polynomials 
$h _m \in k[x,y] $, $m \in \bb Z_{\ge 0}$ 
such that 
\begin{enumerate}
\item $h_{m}$ is divisible by $(x+y)^{m}$,
\item $h_{m}$ lies in $A[[y]]$, 
\item $ I ^{(n +1)} $ is generated by $h_m $, $m \ge n+1$,  
\item  If we write 
\[
 h _m = h_{m,0} + h_{m,1}y +h_{m,2}y^2 + \cdots, \, h _{m,i} \in A,  
\]
then $\ord h_{m,0} = s_m$ and $\ord h_{m,1} \ge s_m -1$.
Moreover $ \ord h_{m,1} = s_m -1 $ if and only if $s_m -1 \in S$.
\end{enumerate}
Being generated by $h_{m,0}$, $m \ge n+1$, the ideal $\f a _{n+1}$ is 
ideantical to the set $\{ f \in A | \ord f \ge s_{n+1}   \}$. 
Since $\c Z_{n,\epsilon}$ corresponds to 
the homomorphism 
$ \f a _{n+1} \to A/\f a_{n+1}$ which takes $h_{m,0}$ to $h_{m,1}$,
 $\c Z_{n,\epsilon }$ is non-trivial if and only if $h_{m,1} \notin \f a _{n+1}$ for some $m \ge n+1$
if and only if $s _{n+1} -1 \in S$.
This completes the proof.
\end{proof}

\begin{definition}
Let $X$ be a reduced scheme and $\nu : Y \to X$
the normalization. Then the {\em conductor ideal}  $\f c \subset \c O_{X}$
is the annihilator of $ \nu_{*}\c O _{Y} / \c O_{X}  $.
The {\em conductor subscheme} of $ X$ is the 
closed subscheme defined by $\f c$.
\end{definition}

By definition, the conductor subscheme is the non-normal locus
endowed with a suitable scheme structure.

\begin{corollary}[\cite{highernash}]\label{cor-curve}
Let $X $ be a  variety of dimension 1, $C$ the conductor subscheme of $X$
 and $ Z \in \Nash _n (X) $ with $Z \not \subset C $. 
Then $\Nash _n (X)$ is normal around $Z$.
\end{corollary}

\begin{proof}[Sketch of the proof]
By \cite[Prop.\ 2.5]{highernash}, we only need to show
the same assertion for formal irreducible curve. 
In this case, with the notations as above,
the conductor ideal $\f c\subset A$ is $\{ f \in A | \ord f \ge c \}$, where
$c := \min \{ s \in S | \forall s' \ge s, \ s' \in S \}$.
Since $\f a_{n+1}=\{f \in A|\ord f \ge s_{n+1}\}$, 
\begin{align*}
\f a _{n+1} \not \supset \f c \Leftrightarrow  
s_{n+1} >c \Rightarrow s_{n+1}-1 \in S.
\end{align*}
By Theorem \ref{thm-formal-curve}, $\Nash _{n} (X)$ is normal.
\end{proof}

Similar arguments apply to the flag higher Nash blowup.

\begin{corollary}\label{cor-formal-curve-flag}
Let the notations be as in Theorem \ref{thm-formal-curve}. 
Then $\fNash _n (\b X)$ is normal if and only if $s_m - 1 \in S$ for some $1 \le m \le n+1$.
\end{corollary}

\begin{proof}
We keep the notation above. The normalization $\b Y \to \b X$ factors as 
$ \b Y \xrightarrow{\psi _n } \fNash _n (\b X) \to \b X$.  It is clear that $\fNash _n (\b X)$
is normal if and only if $\psi_n$ is an isomorphism.
Since $\psi_n$ corresponds to  $( \c Z_0, \dots, \c Z _{n} )$, the last condition is equivalent to that for some $m \le n+1$, 
$\c Z _{m, \epsilon } $  is non-trivial, equivalently, for some $m \le n+1$  $s_m -1 \in S$.
\end{proof}

\begin{corollary}\label{cor-curve-flag}
Let $X$ be a variety of dimension 1, $C $ its conductor subscheme and $ Z _*  \in \fNash _n (X) $ with $Z _{n} \not \supset C $. 
Then $\fNash _n (X)$ is normal around $ Z _*$.
\end{corollary}

\begin{proof}
Arguments similar to the proof of Corollary \ref{cor-curve} apply also to this corollary.
\end{proof}

\subsubsection{Various blowups of $\Spec k[[x^{5},y^{7}]]$}\label{5-7-example}

Let $\b X := \Spec k[[x^{5},x^{7}]]$.
Then
\[
 S = \{0,5,7,10,12,14,15,17,19,20,21,22\}\cup \{ n|n\ge 24\} .
\]
We compute $h_{m}$ ($m \le 6$) in the proof of Theorem \ref{thm-formal-curve} to be
\begin{gather*}
h_{1} = x^{5}+y^{5}, \\
h_{2}= x^7 - (7/5)x^5y^2 - (2/5)y^7, \\
h_{3}= x^{10} + (25/7)x^7y^3 - 3x^5y^5 - (3/7)y^{10}, \\
h_{4}= x^{12 }- (14/5)x^{10}y^2 - 4x^7y^5 + (12/5)x^5y^7 + (1/5)y^{12}, \\
h_{5}=x^{14} - (21/5)x^{12}y^2 + (147/25)x^{10}y^4 + (24/5)x^7y^7 \\ -( 56/25)x^5y^9 - (3/25)y^{14}, \text{ and} \\
h_{6}=x^{15} + (125/49)x^{14}y - (25/7)x^{12}y^3 + 3x^{10}y^5 \\+ (75/49)x^7y^8 - (4/7)x^5y^{10} - (1/49)y^{15} .
\end{gather*}
(Check that for every $i \le 6$, $h_{i} \in k[[x^{5},x^{7},y]]$ and 
$h_{i}$ is divisible by $(x+y)^{i}$.)
If $A_{n}$ is the coordinate ring of $\Nash _{n} (\b X)$,
then $A_{n}$ is (isomorphic to) the least complete $k$-subalgebra
$B \subset k[[y]]$ such that $k[[y^{5},y^{7}]] \subset B$
and for some generators $f_{\lambda}$, $\lambda \in \Lambda $,
 of $I^{(n+1)}$, for every $\lambda \in \Lambda$, $f_{\lambda} \in B$.
If $A'_{n} $ is the coordinate ring of $\fNash _{n} (\b X)$,
then $A'_{n}$ is isomorphic to the least complete $k$-subalgebra
of $k[[y]]$ which contains $A_{m}$, $m \le n$.

Now $I^{(2)}$ is generated by $h_{2}$ and $h_{1}^{2}$.
So $A_{1} = A'_{1} = k[[y^{2},y^{5}]]$.
Then $I^{(3)}$ is generated by $h_{3}$,
$h_{4}+(14/5)y^{2}h_{3}$, and $h_{1}h_{3}$.
So $A_{2} = k[[y^{3},y^{5},y^{7}]]$
and $A'_{2} = k[[y^{2},y^{3}]]$. 
Next $I^{(4)}$ is generated by $h_{4}$, $h_{5}$ and 
\[
h_{6}-(125/49)y h_{5} -(50/7)y^{3}h_{4}
= x^{15} + 8x^{10}y^5 + (125/7)x^7y^8 - 12x^5y^{10} - 8/7y^{15}.
\]
So $A_{3}=k[[y^{2},y^{5}]]$ 
and $A'_{3}=k[[y^{2},y^{3}]]$.
Finally $I^{(5)}$ is generated by $h_{5}$, $h_{6}-(25/49)y h_{5}$ 
and $h_{2}h_{3}$. So $A_{4}=A'_{4}=k[[y^{2},y^{3}]]$.

For $\b X':= \Spec k[[x^{2},x^{5}]]\cong \Nash (\b X)$, 
we similarly define $ I' \subset k[[x^{2},x^{5},y]]$. 
Then $I'^{(2)}$ is generated by
\begin{align*}
  x^{4}-2x^{2}y^{2}+y^{4}  \text{ and }
  x^{5}+(5/2)x^{2}y^{3}-(3/2)y^{5}.
\end{align*}
So if $\Nash^{n}(\cdot)$ denotes the $n$-times iteration of 
the classical Nash blowup, then
 $\Nash^{2} (\b X) = \Nash (\b X') \cong \Spec k[[x^{2},x^{3}]]$,
and $\Nash ^{3}(\b X) \cong \Spec k[[x]]$.

If $B(\cdot)$ denotes the blowup with respect to the
reduced special point and $B^{n}(\cdot)$
is its $n$-times iteration, then it is easy to see that
\begin{gather*}
B(\b X) \cong \Spec k[[x^{2},x^{5}]] \\
 B^{2}(\b X)\cong \Spec k[[x^{2},x^{3}]] \\
B^{3}(\b X) \cong \Spec k[[x]].
\end{gather*}
Thus  the flag $n$-th Nash blowup differs also from
the $n$-times iteration of blowups with respect to 
the reduced special point. 

Every blowup considered here is the spectrum of
the complete algebra associated to a numerical monoid.
Table \ref{table-5-7} shows the correspondence between
 monoids and blowups.

\begin{table}
\begin{center}
\begin{tabular}{|c|c|c|c|c|c|} \hline
n &1&2&3&4&5 \\ \hline
$\Nash _{n} (\b X)$ & $\langle 2,5 \rangle$ & $\langle 3,5,7 \rangle$& $\langle 2,5 \rangle$&$\langle 2,3 \rangle$ &  $\langle 1 \rangle$\\ \hline
$\fNash _{n} (\b X)$ &$\langle 2,5 \rangle$  & $\langle 2,3 \rangle$ &$\langle 2,3 \rangle$  & $\langle 2,3 \rangle$ & $\langle 1 \rangle$ \\ \hline
$\Nash ^{n} (\b X)$ &$\langle 2,5 \rangle$  & $\langle 2,3 \rangle$ &$\langle 1 \rangle$ & $\langle 1 \rangle$&$\langle 1 \rangle$  \\ \hline
$B^{n} (\b X)$ & $\langle 2,5 \rangle$ & $\langle 2,3 \rangle$ &$\langle 1 \rangle$ & $\langle 1 \rangle$&  $\langle 1 \rangle$\\ \hline
\end{tabular}
\caption{The numerical monoids associated to various  blowups of $\Spec k[[x^{5},x^{7}]]$.
Here $\langle a_{1},\dots, a_{l}  \rangle$ is the numerical monoid generated by the natural numbers $a_{1},\dots,a_{l}$.}
\label{table-5-7}
\end{center}
\end{table}

\begin{conjecture}
Let $\b X$ be a formal curve with the associated numerical monoid
$S=\{0= s_{0}<s_{1}< \cdots \}$.
Then the associated numerical monoid of $\Nash _{n} (\b X)$
is generated by  $ s_{m}-s_{l}  $, $ m > n $, $ l \le n$.
That of $\fNash _{n} (\b X)$ is generated by 
$ s_{m}-s_{l}  $, $ m >l  $, $ l \le n$.
\end{conjecture}

\subsection{Product of curves}

Let $\b X_{i} = \Spec A_i$, $i=1,2,\dots, l$, be formal irreducible curves. 
As above, we fix embeddings $A_i \hookrightarrow k[[x_i]]$, define their associated
 numerical monoids $S_i$ and write
\[
 S_{i} = \{ 0 = s_{i,0} < s_{i,1} < s_{i,2} < \cdots  \}.
\]
Similarly we define normalizations $\b Y_{i}:=\Spec k[[y_{i}]] \to \b X_{i}$.
Since formal curves are algebraizable, 
by Theorem \ref{thm-product},  $\hat \prod_{i} \fNash _{n} (\b X_{i}) \cong 
 \fNash _{n} (\hat \prod_{i} \b X_{i}) $.
So there are natural morphisms 
 \[
 \hat \prod_{i} \b Y_{i} \to \hat \prod_{i} \fNash _{n} (\b X_{i}) \cong 
 \fNash _{n} (\hat \prod_{i} \b X_{i})  \to \Nash _{n} (\hat \prod_{i} \b X_{i}) \to \hat \prod_{i} \b X_{i}.
\]
By Corollary \ref{cor-formal-curve-flag}, we obtain:

\begin{theorem}\label{thm-product-curves-flag}
$\fNash _{n} (\hat \prod _{i} \b X _{i})$ is regular if and only if 
for every $i$, there exists $1 \le m_{i} \le n+1 $ such that $s_{i,m_{i}}-1 \in S_{i}$. 
\end{theorem}

Strangely  the same statement holds for the simple higher Nash blowup
if $l \ge 2$.

\begin{theorem}\label{thm-product-curves-simple}
Suppose $l \ge 2$. Then
$\Nash _{n} (\hat \prod _{i} \b X _{i})$ is regular if and only if 
for every $i$, there exists $m_{i} \le n+1 $ with $s_{i,m_{i}}-1 \in S_{i}$. 
\end{theorem}

\begin{proof}
For each $i$, we define 
$I_{i} \subseteq A_{i}[[y_{i}]]$ and $h_{i,m} \in A_{i}[[y_{i}]]$ $(m=0,1,\dots)$ as in the proof of Theorem \ref{thm-formal-curve}.
We denote by  $\psi _{n}$ the natural morphism $ \hat \prod _{i}  \b Y_{i} \to \Nash _{n} (\hat \prod _{i}  \b X_{i})$.
Then $\Nash _{n} (\hat \prod _{i}  \b X_{i})$ is regular if and only if $ \psi_{n}$ is
an isomorphism. 
Consider standard tangent vectors
\[
 \epsilon _{i}: \Spec k[y_{i}]/(y_{i}^{2}) = \Spec k[y_{1}, \dots, y_{l}]/(y_{1},\dots, y_{i-1},y_{i}^{2}, y_{i+1},\dots,y_{l})
 \to \hat \prod _{i} \b  Y_{i} .
\]
Then $\psi_{n}$ is an isomorphism if and only if the
tangent vectors $ \psi _{n} \circ \epsilon_{i} $ are linearly independent.

Let $B := A_{1} \hat \otimes \cdots \hat \otimes A_{l} \subset k[[x_{1},\dots,x_{l}]]$
and $\widetilde{I_{i} ^{ (m_{i} )  }} := I_{i} ^{ (m_{i} ) }B[[y_{1},\dots,y_{l}]]$.
Then the family corresponding to $\psi_{n}$ is defined by the ideal
\[
J^{(n+1)} = \sum_{\sum m_{i} = n+1}  \prod _{i} \widetilde{I_{i} ^{ (m_{i} )  }} \subset 
 B[[y_{1},\dots,y_{l}]],
\]
which is generated by  $\prod_{ i} h _{i,m _{i}} $,   $\sum _{i} m_{i} \ge n+1$.
Set 
\[
 T_{n+1} := \{ (s_{1,m_{1}}, \dots, s_{l,m_{l}}) |  \sum_{i} m_{i} \ge n+1 \} 
 \subseteq S_{1} \times \cdots \times S_{l} \subseteq \bb Z_{\ge 0} ^{l}.
\]
Let $k[[T_{n+1}]] \subset k[[x_{1},\dots,x_{l}]]$
be the complete semigroup algebra (without unit)
associated to $T_{n+1}$.
For $f \in k[[x_{1},\dots,x_{l}]]$,
we denote by $\initial (f)$ its initial form, that is,
the nonzero homogeneous part of the lowest degree.
Then the special point $W_{n} \in \Nash _{n} (\hat \prod _{i} \b  X_{i})$
is defined by the ideal
\begin{equation}\label{ideal-expression}
\f b _{n+1} := J^{(n+1)} B   = \{ f \in B | 
\initial (f) \in k [[ T_{n+1}  ]] \}.
\end{equation}
The Zariski tangent space of $\Nash _{n} (\hat \prod _{i}  \b X_{i})$ at $W_{n}$
is identified with a $k$-subspace of
\[
 \Hom_{B\text{-modules}} ( \f b _{n+1}, B / \f b_{n+1} ).
\]
Writing
\[
 h_{i,m} = h_{i,m,0} + h_{i,m,1} y_{i} + h_{i,m,2} y_{i}^{2} + \cdots, \, h_{i,m,j} \in A_{i},
\]
we have 
\begin{gather*}
\prod_{ i} h _{i,m _{i}} = \prod_{ i} h _{i,m _{i},0} + \sum _{i} \left(  h_{i,m_{i},1}  \prod_{ j \ne i} h _{j,m _{j},0} \right) y_{i} + \cdots, \text{ and} \\
 \f b _{n+1} =  ( \prod_{ i} h _{i,m _{i},0}  | \sum _{i} m_{i} \ge n+1) .
\end{gather*}
So the tangent vector $ \psi_{n} \circ \epsilon _{i} $ is identified with
the homomorphism $ \xi _{i}: \f b _{n+1} \to B / \f b _{n+1}$ that takes
$ \prod_{ j} h _{j,m _{j},0}  $ to $ h_{i,m_{i},1}\prod_{ j \ne i} h _{j,m _{j},0}  $.
From the expression (\ref{ideal-expression}), if $\sum m_{j} >n+1$, we have
$ \xi_{i} (\prod_{ j} h _{j,m _{j},0})  = 0 $.
Moreover when $\sum {m_{j}} =n+1$, 
\begin{align*}
 &\xi_{i} (\prod_{ j} h _{j,m _{j},0})  \ne 0 \\
  &\Leftrightarrow 
 \initial (h_{i,m_{i},1}) = c y_{i}^{s_{m_{i}}-1},\ c \in k \setminus \{0\} \\
 &\Leftrightarrow  s_{i,m_{i}} -1 \in S_{i} .
 \end{align*}

To prove the ``only if'' in the theorem, we now suppose that for some $i$ and every $m \le n+1$,
$s_{i,m} - 1 \notin S_{i}$. Then for any $ ( m_{1}, \dots, m_{l})  $
with $\sum m_{j} = n+1$, we have $\xi _{i} ( \prod_{ j} h _{j,m _{j},0} ) = 0$,
so $\xi_{i}=0$. As a consequence, $\psi _{n}$ is  not an isomorphism
and $\Nash _{n} (\hat \prod X_{i})$ is not regular.

Let us now prove the ``if''. Let $m_{i} \le n+1$, $i=1,\dots,l$, be such that
$s _{m_{i}} -1 \in S_{i}$. Then for each $i$,
there exist $ n_{j} $ $(1 \le j \le l)$ 
such that $n_{i} = m _{i} $
and $ \sum_{j } n_{j}=n+1$. (Here the assumption $l \ge 2$ is necessary.)
We have that
\begin{align*}
 \xi _{j'} ( \prod_{ j} h _{j,n _{j},0} )  = \text{the class of } 
  h _{j,n_{j'},1} \prod_{ j \ne j'} h _{j,n _{j},0} \text{ modulo $\f b_{n+1}$} ,
\end{align*}
which is nonzero if $j' =i$. 
It follows that $\xi _{i} (\prod_{ j} h _{j,n _{j},0} )$
is not any $k$-linear combination of $ \xi _{j'} (\prod_{ j} h _{j,n _{j},0} ) $, $j' \ne i$. 
So
if a linear relation $ \sum  _{j} r_{j} \xi _{j}$ ($r_{j} \in k$) holds,
then $r _{i} = 0$. Since this holds for every $i$, the $\xi_{i}$'s are
linearly independent, which completes the proof.
\end{proof}

\begin{example}\label{expl-not-isom}
Suppose that $l=2$, $S_{1} = \bb Z _{\ge 0}$ and 
\[
S_{2} := \langle 3,4 \rangle = \{0,3,4,6,7,8,9,\dots\}.
\]
Then $ \Nash _{n} ( \b X_{1})\hat \times \Nash_{n} ( \b X_{2})$ is regular
if and only if $\Nash_{n} ( \b X_{2})$ is regular 
if and only if $n \ne 0,2$. On the other hand, 
$ \Nash _{n} ( \b X_{1} \hat \times  \b X_{2})  $ is regular if and only if
$n > 0 $. So 
\[
\Nash _{2} ( \b X_{1})\hat \times \Nash_{2} ( \b X_{2})
\not \cong   \Nash _{2} ( \b X_{1} \hat \times  \b X_{2})  .
\]
In particular, this example says that the simple higher Nash blowup
is not generally compatible with either products or smooth morphisms.
\end{example}

\begin{corollary}
Let $X$ be a variety such that 
every analytic branch of $X$ at every point has the same singularity type
as the product of curves, and  $C \subset X$ the conductor subscheme.
Then for every $n$, the normalization $Y \to X$
factors as $Y \to \fNash _{n}(X) \to \Nash _{n}(X) \to X$.
Moreover for $Z \in \Nash _{n} (X)$ 
with $Z \not \subset C$, (resp.\ $Z _{*} \in \fNash _{n}$
with  $Z _{n} \not \subset C$),
 $\Nash _{n} (X)$ (resp.\ $\fNash _{n} (X)$) is smooth
at $Z$ (resp.\ $Z_{*}$).
 \end{corollary}

\begin{proof}
From \cite[Prop.\ 2.5]{highernash}, it suffices to show 
the same assertion for the product of formal irreducible curves. 
Let the notations be as in the proof of Theorem \ref{thm-product-curves-simple}.
We may suppose that $l \ge 2$.
Then the special points of $ \Nash _{n} (\hat \prod_{i} \b X_{i}) $ and $\fNash _{n} (\hat \prod_{i} \b X_{i})$
correspond to the same ideal $\f a _{n +1}$.
Let  $ \f c _{i}  $, $i=1,\dots,l$, be the conductor
ideals of $A_{i}$, and $\tilde{\f c_{i}} := \f c_{i} B$.
 Then the conductor ideal
of $B$ is $\f c := \prod \tilde{ \f c _{i}}$.
Let $U := \{ (m_{1},\dots,m_{l}) | m_{i} \ge c_{i} \} \subset \bb Z_{\ge 0}^{l}$.
Then $\f c=k[[U]] \subset k[[x_{1},\dots,x_{l}]]$.

Now the assumption $ \f a _{n+1} \not \supset \f c $ is equivalent to
that $T_{n+1} \not \supset U$. If it is the case, then  for every $i$, 
$c _{i}+1 \le   s _{i,n+1}$. From Theorems \ref{thm-product-curves-flag}
and \ref{thm-product-curves-simple}, $\Nash _{n} (\hat \prod _{i}\b X_{i})$ and
$\fNash _{n} (\hat \prod_{i} \b X_{i})$ are regular.
\end{proof}

Lastly we raise two problems.

\begin{problem}
Is the correspondence of
Theorems \ref{thm-product-curves-flag} and \ref{thm-product-curves-simple}
only a coincidence?
Or, could it be that $\fNash _{n} (X) \cong \Nash _{n} (X)$ 
for any variety $X$ of dimension $\ge 2$?
\end{problem}

\begin{problem}
Let $X$ be a variety,  $C$ its conductor subscheme,
$Z \in \Nash _{n} (X)$ with $Z \not \subset C$ and $Z_{*} \in \fNash _{n}(X)$
with $Z_{n} \not \subset C$.
Then  are $\fNash _{n} (X) $ and $\fNash _{n}(X)$ normal at $Z$
and $Z_{*}$ respectively?
\end{problem}


\begin{thebibliography}{9}


\bibitem{BierstoneMilmanCanonicalResolution}
E.\ Bierstone and P.D.\ Milman.
\newblock Canonical desingularization in characteristic zero by  blowing up the maximal strata
of a local invariant.
\newblock {\em Invent.\ Math.}, 128:207--302, 1997.


\bibitem{EncinasVillamayorActa}
S.\ Encinas and O.\ Villamayor.
\newblock {Good points and constructive resolution of singularities}.
\newblock {\em Acta Math.}, 181(1):109--158, 1998.


\bibitem{VillamayorConstructiveness}
O.\ Villamayor.
\newblock Constructiveness of Hironaka's resolution.
\newblock {\em Ann.\ Sci.\ \'Ecole Norm.\ Sup.\ (4)}, 22: 1--32, 1989.

\bibitem{VillamayorPatching}
O.\ Villamayor.
\newblock Patching local uniformizations.
\newblock {\em Ann.\ Sci.\ \'Ecole Norm.\ Sup.\ (4)}, 25: 629--677, 1992.


\bibitem{highernash} T.\ Yasuda, Higher Nash blowups, math.AG/0512184.
\end{thebibliography}
\end{document}